\documentclass{commat}

\usepackage{tikz}
\usetikzlibrary{arrows,calc,positioning}
\tikzstyle{intt}=[draw,text centered,minimum size=6em,text width=5.25cm,text height=0.34cm]
\tikzstyle{intl}=[draw,text centered,minimum size=2em,text width=4cm,text height=0.34cm]
\tikzstyle{int}=[draw,minimum size=2.5em,text centered,text width=5cm]
\tikzstyle{intg}=[draw,minimum size=3em,text centered,text width=2.5cm]
\tikzstyle{intgg}=[draw,minimum size=3.6em,text centered,text width=6.cm]
\tikzstyle{sum}=[draw,shape=circle,inner sep=2pt,text centered,node distance=3.5cm]
\tikzstyle{summ}=[drawshape=circle,inner sep=4pt,text centered,node distance=3.cm]
\tikzstyle{intr}=[draw,minimum size=6.5em,text centered,text width=4.5cm]
\tikzstyle{intrr}=[draw,minimum size=2.5em,text centered,text width=4.5cm]
\tikzstyle{intrrr}=[draw,minimum size=5.5em,text centered,text width=4.5cm]
\tikzstyle{intaaa}=[draw,minimum size=2em,text centered,text width=3.5cm]
\tikzstyle{inta}=[draw,minimum size=1.5em,text centered,text width=3cm]
\tikzstyle{intaa}=[draw,minimum size=1.5em,text centered,text width=2.5cm]

\title{%
    On the matrix function $_pR_q(A, B; z)$ and its fractional calculus properties    
    }

\author{%
    Ravi Dwivedi and Reshma Sanjhira
    }

\affiliation{
    \address{Ravi Dwivedi --
    Department of Science, Govt Naveen College Bhairamgarh, Bijapur (CG), 494450, India.
    }
    \email{%
    dwivedir999@gmail.com
    }
    \address{Reshma Sanjhira --
    Department of Mathematics, Faculty of Science, The Maharaja Sayajirao University of Baroda, Vadodara-390002, India. Department of Mathematical Sciences, P. D. Patel Institute of Applied Sciences, Faculty of Applied Sciences, Charotar University of Science and Technology, Changa-388 421, India.
        }
    \email{%
    reshmasanjhira1508@gmail.com
    }
    }

\abstract{%
    The main objective of the present paper is to introduce and study the function $_pR_q(A, B; z)$ with matrix parameters and  investigate the convergence of this matrix function. The contiguous matrix function relations, differential formulas and the integral representation for the matrix function $_pR_q(A, B; z)$ are derived. Certain properties of the matrix function $_pR_q(A, B; z)$  have also been studied from fractional calculus  point of view. Finally, we emphasize on the special cases namely the generalized matrix $M$-series, the Mittag-Leffler matrix function and its  generalizations and some matrix polynomials.
    }



\VOLUME{31}
\YEAR{2023}
\NUMBER{1}
\firstpage{43}
\DOI{https://doi.org/10.46298/cm.10205}

\begin{paper}

\section{Introduction}\label{s:i}

Special matrix functions play an important role in mathematics and physics. In particular, special matrix functions appear in the study of statistics \cite{cm}, probability theory \cite{js} and Lie theory \cite{ds6}, \cite{atj}, to name a few. The theory of special matrix functions has been initiated by J\'odar and Cort\'es who studied matrix analogues of gamma, beta and Gauss hypergeometric functions \cite{jjc98a}, \cite{jjc98b}. Dwivedi and Sahai generalized the study of one variable special matrix functions to $n$-variables \cite{ds1}-\cite{ds5}. Some of the extended work of Appell matrix functions have been given in \cite{al}. Certain polynomials in one or more variables have been introduced and studied from matrix point of view, see \cite{acc}, \cite{ars}, \cite{ca}, \cite{djl}, \cite{sb}, \cite{sd06}. Recently, the  generalized Mittag-Leffler matrix function have been introduced and studied in \cite{snd}.

It appears from the literature that the function $_pR_q(\alpha, \beta; z)$ were systematically studied in \cite {ds}. In this article, we introduce a new class of matrix function, namely $_pR_q(A, B; z)$ and discuss its regions of convergence. We also give  contiguous matrix function relations, integral representations and differential formulas satisfied by the matrix function $_pR_q(A, B; z)$. The matrix analogues of generalized $M$-series $_pM_q^{\alpha, \beta}(\gamma_1, \dots, \gamma_p,\delta_1, \dots, \delta_q; z)$, Mittag-Leffler functions and its generalizations have been presented as special cases of the matrix function $_pR_q(A, B; z)$. 
The paper is organized as follows:

In Section~2, we list the basic definitions and results from special matrix functions that are needed in the sequel. In Section~3, we introduce the  matrix function $_pR_q(A, B; z)$ and prove a theorem on its absolute convergence. In Section~4, we give contiguous matrix function relations and differential formulas satisfied by $_pR_q(A, B; z)$. In Section~5, an integral representation of the matrix function $_pR_q(A, B; z)$ motivated by the integral of beta matrix function has been given. In Section~6, the fractional order integral and differential transforms of the matrix function  $_pR_q(A, B; z)$ have been determined. Finally, in Section~7, we present the Gauss hypergeometric matrix function and its generalization, the matrix $M$-series, the Mittag-Leffler matrix function and its generalizations  and some matrix polynomials as special cases of $_pR_q(A, B; z)$. 
	
\section{Preliminaries} 
Let the spectrum of a matrix $A$ in $\mathbb{C}^{r\times r}$, denoted by $\sigma(A)$, be the set of all eigenvalues of $A$.
Recall that a matrix $A \in \mathbb{C}^{r \times r}$ is said to be positive stable when
\[
    \beta(A)
    = \min \{\,\Re(z) \mid z \in \sigma(A)\,\} > 0.
\]
For a positive stable matrix $A \in \mathbb{C}^{r \times r}$, the gamma matrix function is defined by \cite{jjc98a}
\[ \Gamma(A) = \int_{0}^{\infty} e^{-t} \, t^{A-I}\, dt
\]
and the reciprocal gamma matrix function is defined as \cite{jjc98a}
	  \begin{equation}
	 \Gamma^{-1}(A)= A(A+I)\dots (A+(n-1)I)\Gamma^{-1}(A+nI) , \  n\geq 1.\label{eq.07}
	  \end{equation}
The Pochhammer symbol  for  $A\in \mathbb{C}^{r\times r}$ is given by \cite{jjc98b}
	  \begin{equation}
	  (A)_n = \begin{cases}
	  I, & \text{if $n = 0$,}\\
	  A(A+I) \dots (A+(n-1)I), & \text{if $n\geq 1$}.
	  \end{cases}\label{c1eq.09}
	  \end{equation}
	  This gives
	  	   \begin{equation}
	   (A)_n = \Gamma^{-1}(A) \ \Gamma (A+nI), \qquad n\geq 1.\label{c1eq.010}
	   \end{equation} 
If $A \in \mathbb{C}^{r\times r}$ is a positive stable matrix  and $n\geq 1$ is an integer, then the  gamma matrix function can also be defined in the form of a limit as \cite{jjc98a}
\begin{equation} 
\Gamma (A) = \lim_{n \to \infty} (n-1)! \, (A)_n^{-1} \, n^A. \label{eq10} 
\end{equation}
If $A$ and $B$ are positive stable matrices in $\mathbb{C}^{r \times r}$, then the beta matrix function is defined as \cite{jjc98a}
\begin{equation}
\mathfrak{B}(A,B) =\int_{0}^{1} t^{A-I} \, (1-t)^{B-I} dt.\label{c1eq11}
\end{equation}
Furthermore, if $A$, $B$ and $A+B$ are positive stable matrices in $\mathbb{C}^{r \times r}$ such that $AB = BA$, then the beta matrix function is defined as \cite{jjc98a} 
\begin{equation}
\mathfrak{B}(A,B) = \Gamma(A)\,\Gamma(B)\,\Gamma^{-1} (A+B).
\end{equation}
	 Using the Schur decomposition of $A$,  it follows that \cite{gl}, \cite{vl}
	  \begin{equation}
	  \Vert e^{tA}\Vert \leq e^{t\alpha(A)} \sum_{k=0}^{r-1}\frac{(\Vert A\Vert r^{1/2} t)^k}{k!}, \ \ t\geq 0.\label{eq09}
	  \end{equation}
	  We shall use the notation $\Gamma \left(\begin{array}{c}
	 A_1, \dots, A_p \\
	 B_1, \dots, B_q
	  \end{array}\right)$ for  $\Gamma (A_1) \cdots \Gamma (A_p) \, \Gamma ^{-1}(B_1) \cdots \Gamma ^{-1} (B_q)$.
			
\section{The matrix function $_pR_q(A, B; z)$}	
J\'odar and Cort\'es \cite{jjc98b} defined the Gauss hypergeometric function with matrix parameters denoted by $_2F_1(A,B;C;z)$, where $A$, $B$, $C$ are matrices in $\mathbb{C}^{r\times r}$, and determined its region of convergence and integral representation. A natural generalization of the Gauss hypergeometric matrix function is obtained in \cite{ds1} by introducing an arbitrary number of matrices as parameters in the numerator and denominator and referring to this generalization as the generalized hypergeometric matrix function, $_pF_q(A_1,\dots,A_p;B_1,\dots,B_q;z)$. We now give an extension of the generalized hypergeometric matrix function. Let $A$, $B$, $C_i$ and $D_j$, $1\leq i\leq p$, $1\leq j\leq q$, be matrices  in $\mathbb{C}^{r\times r}$ such that $D_j + kI$ are invertible for all integers $k\geq 0$. Then, we define the matrix function $_pR_q(A, B; z)$ as
	\begin{align}
_pR_q(A, B; z) &= {}_pR_q\left(\begin{array}{c}
	C_1, \dots, C_p \\
	D_1, \dots, D_q
	\end{array}\mid A, B; z\right)\nonumber\\
	 & = \sum_{n\geq 0} \Gamma^{-1}(nA+B) \, (C_1)_n\dots (C_p)_n \, (D_1)_n^{-1}\dots (D_q)_n^{-1} \, \frac{z^n}{n!}.\label{3.1}
	\end{align}
In the following theorem, we find the regions in which the matrix function $_pR_q(A, B; z)$ either converges or diverges.
	 \begin{theorem}
	 	Let $A, B, C_1,\dots,C_p, D_1,\dots ,D_q$ be positive stable matrices in $\mathbb{C}^{r\times r}$. Then the matrix function $_pR_q(A, B; z)$ defined in \eqref{3.1} converges or diverges in one of the following regions:
	 	\begin{enumerate}
	 	\item If $p \le q+1$, the matrix function converges absolutely for all finite $z$.
	 	\item If $p = q + 2$, function converges for $\vert z \vert < 1$ and diverges for $\vert z \vert > 1$.
	 	\item If $p = q + 2$ and $\vert z \vert = 1$, the function converges absolutely for
        \[ \beta (D_1) + \cdots + \beta(D_q) > \alpha (C_1) + \cdots + \alpha (C_p). \]
	 	\item If $p > q + 2$, the function diverges for all $z \ne 0$.
	 	\end{enumerate}
	 \end{theorem}
	 \begin{proof} 
	 	Let $U_n(z)$ denote the general term of the series \eqref{3.1}. Then, we have
	 	\begin{align}
	 	\Vert U_n(z)\Vert &\le \Vert \Gamma^{-1} (nA+B)\Vert \, \prod_{i =1}^{p} \Vert (C_i)_n\Vert \, \prod_{j =1}^{q} \Vert (D_j)_n^{-1}\Vert \, \frac{\vert z\vert^n}{n!}\nonumber\\
	 	& \le \Vert \Gamma^{-1} (nA+B)\Vert \, \prod_{i =1}^{p} \left\Vert \frac{(C_i)_n n^{C_i} n^{-C_i} (n-1)!}{(n-1)!} \right\Vert \, \nonumber\\
	 	& \quad \times \prod_{j =1}^{q} \left\Vert \frac{(D_j)_n^{-1} n^{D_j} n^{-D_j} (n-1)!}{(n-1)!} \right\Vert \, \frac{\vert z\vert^n}{n!}.
	 	\end{align}
The limit definition of gamma matrix function \eqref{eq10} and Schur decomposition \eqref{eq09} yield
\begin{align}
\Vert U_n(z)\Vert & \le N \ S \ ((n-1)!)^{p-q-2} \ n^{\sum_{i= 1}^{p} \alpha(C_i) - \sum_{j= 1}^{q} \beta(D_j) - 1} \, \vert z\vert^n, \label{3.3}
\end{align}
where $N = \Vert \Gamma^{-1} (C_1)\Vert  \cdots \Vert\Gamma^{-1}(C_p)\Vert \Vert\Gamma(D_1)\Vert \cdots \Vert\Gamma(D_q)\Vert$ and 
\begin{align}
S = \left(\sum_{k= 0}^{r-1} \frac{(\max \{\Vert C_1\Vert, \dots, \Vert C_p\Vert, \Vert D_1\Vert, \dots, \Vert D_q\Vert\} \, r^{\frac{1}{2}} \, \ln n)^k}{k !}\right)^{p+q}.
\end{align}
Thus, it can be easily calculated from \eqref{3.3} and comparison theorem of numerical series that the matrix series \eqref{3.1} converges or diverges in one of the region listed in Theorem~3.1. 
\end{proof}
\section{Contiguous matrix function relations}
In this section, we shall obtain contiguous matrix function relations and differential formulas satisfied by the matrix function $_pR_q(A, B; z)$. The following abbreviated notations will be used throughout the subsequent sections:
\begin{align}
 &R = {}_pR_q(A, B; z) = {}_pR_q\left(\begin{array}{c}
 	C_1, \dots, C_p \\
 	D_1, \dots, D_q
 \end{array}\mid A, B; z\right),\nonumber\\
&R(C_i +) = {}_pR_q\left(\begin{array}{c}
	C_1, \dots, C_{i-1}, C_i + I, C_{i+1}, \dots, C_p \\
	D_1, \dots, D_q
\end{array}\mid A, B; z\right),\nonumber\\
&R(C_i -) = {}_pR_q\left(\begin{array}{c}
C_1, \dots, C_{i-1}, C_i - I, C_{i+1}, \dots, C_p \\
D_1, \dots, D_q
\end{array}\mid A, B; z\right),\nonumber\\
&R(D_j -) = {}_pR_q\left(\begin{array}{c}
C_1, \dots, C_p \\
D_1, \dots, D_{j-1}, D_j - I, D_{j+1}, \dots, D_q
\end{array}\mid A, B; z\right),\nonumber\\
 &{}_{p}R_{q}(A, B+I; z)  = {}_pR_q\left(\begin{array}{c}
C_1, \dots, C_p \\
D_1, \dots, D_q
\end{array}\mid A, B+I; z\right),\nonumber\\
& {}_{p}R_{q}(A, B-I; z)  = {}_pR_q\left(\begin{array}{c}
C_1, \dots, C_p \\
D_1, \dots, D_q
\end{array}\mid A, B-I; z\right).
\end{align}  
Following Desai and Shukla \cite{ds}, we can find $(p+q-1)$ contiguous matrix function relations of bilateral type that connect either $R$, $R(C_1+)$ and $R(C_i+)$, $1 \le i \le p$ or $R$, $R(C_1+)$ and $R(D_j-)$, $1 \le j \le q$. Let $C_i, 1 \le i \le p$ be positive stable matrices in $\mathbb{C}^{r \times r}$ such that $C_i C_k = C_k C_i, 1 \le k\le p, k < i$, $C_i A = AC_i$ and $C_i B = B C_i$. Then, we have 
\begin{align}
R(C_i+) = \sum_{n\geq 0} C_i^{-1} (C_i + nI)  \Gamma^{-1}(nA+B) \, (C_1)_n\cdots (C_p)_n \, (D_1)_n^{-1}\cdots (D_q)_n^{-1} \, \frac{z^n}{n!}.\label{3.6}
\end{align}
If $\theta = z\frac{d}{dz}$ is a differential operator, then we get
\begin{align}
(\theta + C_i) R = \sum_{n\geq 1} (C_i + nI)  \Gamma^{-1}(nA+B) \, (C_1)_n\cdots (C_p)_n \, (D_1)_n^{-1}\cdots (D_q)_n^{-1} \, \frac{z^n}{n!}.\label{3.7}
\end{align}
Equations \eqref{3.6} and \eqref{3.7} together yield
\begin{equation}
(\theta + C_i) R = C_i \, R(C_i+), \quad i = 1, \dots, p.\label{3.8}
\end{equation}
In particular, for $i = 1$, we write
\begin{equation}
(\theta + C_1) R = C_1 \, R(C_1+).\label{3.9}
\end{equation}
Similarly for matrices $D_j \in \mathbb{C}^{r \times r}, 1 \le j \le q$ such that $D_j D_k = D_k D_j, 1 \le k \le q,  k > j$, we obtain a set of $q$ equations, given by
\begin{equation}
\theta \, R + R \, (D_j - I) = R(D_j -) (D_j - I). \label{3.10}
\end{equation}
Now, eliminating $\theta$ from \eqref{3.8} and \eqref{3.10} gives rise to $(p + q - 1)$ contiguous matrix function relations of bilateral type
\begin{equation}
 C_i \, R - R \, (D_j - I) = C_i \, R(C_i+) - R(D_j -) (D_j - I), \ 1 \le i \le p, 1 \le j \le q.
\end{equation}
Equations \eqref{3.8} and \eqref{3.9} produce $(p-1)$ contiguous matrix function relations 
\begin{equation}
(C_1 - C_i) R = C_1 R(C_1 +) - C_i \, R(C_i+), \quad i = 2, \dots, p.\label{4.8}
\end{equation}
Furthermore, Equations \eqref{3.9} and \eqref{3.10} leads to $q$ contiguous matrix function relations 
\begin{equation}
C_1 \, R - R \, (D_j - I) = C_1 \, R(C_1+) - R(D_j -) (D_j - I), \ 1 \le j \le q.\label{4.9}
\end{equation}
The set of matrix function relations given in \eqref{4.8} and \eqref{4.9} are simple contiguous matrix function relations.

Next, we give  matrix differential formulas  satisfied by the matrix function ${}_{p}R_{q}(A, B; z)$.

\subsection{Matrix differential formulas}
\begin{theorem}\label{et1}
Let $A$, $B$, $C_1$, $\dots, C_p$, $D_1$, $\dots, D_q \in \mathbb{C}^{r \times r}$ such that each $D_j + kI, 1\le j \le q$ is invertible for all integers $k\ge 0$. Then the matrix function ${}_{p}R_{q}(A, B; z)$ satisfies the matrix differential formulas
\end{theorem}
\begin{align}
\left(\frac{d}{dz}\right)^r {}_{p}R_{q}(A, B; z) & = (C_1)_r \cdots (C_p)_r \ _pR_q\left(\begin{array}{c}
C_1 + rI, \dots, C_p +rI \\
D_1+rI, \dots, D_q+rI
\end{array}\mid A, rA+B; z\right) \nonumber\\
&  \quad \times  (D_1)_r^{-1} \cdots (D_q)_r^{-1}, \, C_l C_m = C_m C_l, \, C_l A = A C_l, \, C_l B = B C_l, \nonumber\\
&\qquad D_i D_j = D_j D_i, 1 \le l, m \le p,  \, 1 \le i, j \le q;\label{4.10}
\end{align}
\begin{align} 
\left(\frac{d}{dz}\right)^r ({}_{p}R_{q}(A, B; z) z^{D_j-I})
&  = {}_pR_q\left(\begin{array}{c}
C_1, \dots, C_p \\
D_1, \dots, D_{j-1}, D_j - rI, D_{j+1}, \dots, D_q
\end{array}\mid A, B; z\right)\nonumber\\
& \quad \times (-1)^r z^{D_j -(r+1)I} (I-D_j)_r, \  D_i D_j = D_j D_i;\label{4.11}
\end{align}
\begin{align} 
& \left(z^2 \,\frac{d}{dz}\right)^r (z^{C_i-(r-1)I} {}_{p}R_{q}(A, B; z) )\nonumber\\
& = (C_i)_r \, z^{C_i + rI} {}_pR_q\left(\begin{array}{c}
C_1, \dots, C_{i-1}, C_i + rI, C_{i+1}, \dots, C_p \\
D_1, \dots, D_q
\end{array}\mid A, B; z\right), C_i C_j = C_j C_i\nonumber\\
& \qquad C_i A = A C_i, \, C_i B = B C_i, \, 1 \le i, j \le p.\label{4.12}
\end{align}
\begin{proof}
	Differentiating the Equation \eqref{3.1} with respect to $z$, we get
	\begin{align}
\frac{d}{dz} \  _{p}R_{q}(A, B; z) &= \sum_{n\geq 1} \Gamma^{-1}(nA+B) \, (C_1)_n\dots (C_p)_n \, (D_1)_n^{-1}\dots (D_q)_n^{-1} \, \frac{z^{n-1}}{(n-1)!}\nonumber\\
& = \sum_{n\geq 0} \Gamma^{-1}(nA+A+B) \, (C_1)_{n+1}\dots (C_p)_{n+1} \, (D_1)_{n+1}^{-1}\dots (D_q)_{n+1}^{-1} \, \frac{z^{n}}{n!}\nonumber\\
& = (C_1)_1 \cdots (C_p)_1 \ _pR_q\left(\begin{array}{c}
C_1 + I, \dots, C_p +I \\
D_1+I, \dots, D_q+I
\end{array}\mid A, A+B; z\right) \nonumber\\
&  \quad \times  (D_1)_1^{-1} \cdots (D_q)_1^{-1}.
	\end{align}
Proceeding similarly $r$-times, we get the required relation \eqref{4.10}. Using the commutativity of matrices considered in the hypothesis and the way \eqref{4.10} is proved, we are able to prove \eqref{4.11} and \eqref{4.12}.
\end{proof}
\begin{theorem}
	Let $A$, $B$, $C_1$, $\dots, C_p$, $D_1$, $\dots, D_q \in \mathbb{C}^{r \times r}$ such that each $D_j + kI, 1\le j \le q$ is invertible for all integers $k\ge 0$  and $A$, $B-I$ are positive stable. Then the matrix function ${}_{p}R_{q}(A, B; z)$ defined in \eqref{3.1} satisfies the matrix differential formula
	\begin{align}
	zA\frac{d}{dz} \, {}_{p}R_{q}(A, B; z) = {}_{p}R_{q}(A, B-I; z) - (B-I) \, {}_{p}R_{q}(A, B; z), \quad AB = BA.\label{3.23}
	\end{align} 
\end{theorem} 
\begin{proof}
	Using the definition of matrix function ${}_{p}R_{q}(A, B; z)$ and $z\frac{d}{dz} z^n = n z^n$ in the left hand side of \eqref{3.23}, we get 
	\begin{align}
	zA\frac{d}{dz} \, {}_{p}R_{q}(A, B; z) & = \sum_{n\geq 0} nA \Gamma^{-1}(nA+B) \, (C_1)_n\dots (C_p)_n \, (D_1)_n^{-1}\dots (D_q)_n^{-1} \, \frac{z^n}{n!}\nonumber\\
	& = \sum_{n\geq 0} \Gamma^{-1}(nA+B - I) \, (C_1)_n\dots (C_p)_n \, (D_1)_n^{-1}\dots (D_q)_n^{-1} \, \frac{z^n}{n!}\nonumber\\
	& \quad - (B-I) \sum_{n\geq 0} \Gamma^{-1}(nA+B) \, (C_1)_n\dots (C_p)_n \, (D_1)_n^{-1}\dots (D_q)_n^{-1} \nonumber\\
	& \quad \times  \frac{z^n}{n!}, \quad AB = BA\nonumber\\
	& = {}_{p}R_{q}(A, B-I; z) - (B-I) \, {}_{p}R_{q}(A, B; z).
	\end{align}
This completes the proof of \eqref{3.23}.
\end{proof}

\section{Integral representation}
We now find an integral representation of the matrix function ${}_pR_q(A, B; z)$ using the integral of the beta matrix function.
\begin{theorem}
	Let $A$, $B$, $C_1$, $\dots, C_p$, $D_1, \dots, D_q$ be matrices in $\mathbb{C}^{r \times r}$ such that: $C_p$, $D_q$, $D_q - C_p$ are positive stable and $C_p D_j = D_j C_p$ for all $1 \le j \le q$. Then, for $\vert z\vert < 1$, the matrix function ${}_pR_q(A, B; z)$ defined in \eqref{3.1} can be presented in integral form as
	\begin{align}
{}_{p}R_{q}(A, B; z) &= \int_{0}^{1} {}_{p-1}R_{q-1}\left(\begin{array}{c}
C_1, \dots, C_{p-1} \\
D_1, \dots, D_{q-1}
\end{array}\mid A, B; t\,z\right) t^{C_p-I} (1-t)^{D_q-C_p-I} dt\nonumber\\
& \quad \times \Gamma\left(\begin{array}{c}
D_q \\
C_p, D_q-C_p
\end{array}\right).\label{3.14}
	\end{align} 
\end{theorem}
 \begin{proof}
  Since $C_p, D_q, D_q - C_p$ are positive stable and $C_p D_q = D_q C_p$, we have \cite{jjc98b}
  \begin{align}
  (C_p)_n (D_q)_n^{-1} = \left(\int_{0}^{1} t^{C_p + (n-1)I} (1-t)^{D_q-C_p-I} dt \right)  \Gamma\left(\begin{array}{c}
  D_q \\
  C_p, D_q-C_p
  \end{array}\right).\label{3.15}
  \end{align}
  Using \eqref{3.15} in \eqref{3.1}, we get
  \begin{align}
  {}_{p}R_{q}(A, B; z) &=  \sum_{n\geq 0} \int_{0}^{1} \Gamma^{-1} (nA+B) \, (C_1)_n \cdots (C_{p-1})_n \, (D_1)_n^{-1}\cdots (D_{q-1})_n^{-1} \nonumber\\
  &  \quad \times \frac{z^n}{n!} \, t^{C_p + (n-1)I} \, (1-t)^{D_q-C_p-I} dt \ \Gamma\left(\begin{array}{c}
  D_q \\
  C_p, D_q-C_p
  \end{array}\right).\label{3.16}
  \end{align}
  To interchange the integral and summation, consider the product of matrix functions
  \begin{align}
  S_n(z, t) &=  \Gamma^{-1} (nA+B) \, (C_1)_n \cdots (C_{p-1})_n \, (D_1)_n^{-1}\cdots (D_{q-1})_n^{-1} \, \frac{z^n}{n!} \, t^{C_p + (n-1)I}\nonumber\\
  &  \quad \times (1-t)^{D_q-C_p-I} \ \Gamma\left(\begin{array}{c}
  D_q \\
  C_p, D_q-C_p
  \end{array}\right).
  \end{align}
  For $0 < t < 1$ and $n \ge 0$, we get
   \begin{align}
  &\Vert S_n(z, t)\Vert\nonumber\\
   &\le  \left \Vert \Gamma\left(\begin{array}{c}
  D_q \\
  C_p, D_q-C_p
  \end{array}\right)\right \Vert \, \left\Vert\Gamma^{-1} (nA+B) \, (C_1)_n \cdots (C_{p-1})_n \, (D_1)_n^{-1}\cdots (D_{q-1})_n^{-1} \, \frac{z^n}{n!}\right\Vert \, \nonumber\\
  &  \quad \times \Vert t^{C_p -I}\Vert \Vert(1-t)^{D_q-C_p-I}\Vert.
  \end{align}
  The Schur decomposition \eqref{eq09} yields
\begin{align}
\Vert t^{C_p-I}\Vert \ \Vert (1-t)^{D_q-C_p-I}\Vert
&  \leq t^{\alpha(C_p)-1}(1-t)^{{\alpha(D_q - C_p)}-1} \left(\sum_{k=0}^{r-1}
\frac{(\Vert C_p-I\Vert \ r^{1/2} \ \ln t )^k}{k!}\right)\nonumber\\
&  \quad \times \left(\sum_{k=0}^{r-1} \frac{(\Vert {D_q-C_p-I}\Vert \ r^{1/2} \ \ln {(1-t)} )^k}{k!}\right).
\end{align}
Since $0 < t < 1$, we have
\begin{align}
\Vert t^{C_p-I}\Vert \ \Vert (1-t)^{D_q-C_p-I}\Vert \leq \mathcal{A} \  t^{\alpha(C_p)-1}(1-t)^{{\alpha(D_q - C_p)}-1}, 
\end{align}
where
\begin{align}
\mathcal{A} = \left(\sum_{k=0}^{r-1} \frac{(\max \{\Vert C_p-I\Vert, \Vert D_q-C_p - I\Vert\} \ r^{1/2} )^k}{k!}\right)^2.
\end{align}
The matrix series $\Gamma^{-1} (nA+B) \, (C_1)_n \cdots (C_{p-1})_n \, (D_1)_n^{-1}\cdots (D_{q-1})_n^{-1} \, \frac{z^n}{n!}$ converges absolutely for $p \le q + 2$ and $\vert z \vert < 1$; suppose it converges to $S'$. Thus, we get
\begin{align}
\sum_{n\ge 0} \Vert S_n(z, t)\Vert \le f(t) = N S' \mathcal{A} \, t^{\alpha(C_p)-1} \, (1-t)^{\alpha(D_q - C_p)-1}. 
\end{align} 
Since $\alpha(C_p), \alpha(D_q - C_p) > 0$, the function $f(t)$ is integrable and by the dominated convergence theorem \cite{gf}, the summation and the integral can be interchanged in \eqref{3.16}. Using $C_p D_j = D_j C_p, \ 1 \le j \le q$, we get \eqref{3.14}.
\end{proof}
 
\section{Fractional calculus of the matrix function ${}_{p}R_{q}(A, B; z)$}
Let $x > 0$ and $\mu \in \mathbb{C}$ such that $\Re(\mu) > 0$. Then the Riemann-Liouville type fractional order integral and derivatives of order $\mu$ are given by \cite{kst}, \cite{skm}
\begin{equation}
({\bf I}^\mu_a f)(x) = \frac{1}{\Gamma(\mu)} \int_{a}^{x} (x-t)^{\mu - 1} f(t) dt\label{6.1}
\end{equation}
and 
\begin{equation}
{\bf D}_a^\mu f(x) = ({\bf I}^{n-\mu}_a {\bf D}^n f(x)), \quad {\bf D} = \frac{d}{dx}.\label{6.2}
\end{equation}
Bakhet and his co-workers, \cite{ab},  studied the fractional order integrals and derivatives of  Wright hypergeometric and incomplete Wright hypergeometric matrix functions using the operators \eqref{6.1} and \eqref{6.2}. To obtain such they used the following lemma:
\begin{lemma}\label{l1}
Let $A$ be a positive stable matrix in $\mathbb{C}^{r \times r}$ and $\mu \in \mathbb{C}$ such that $\Re(\mu) > 0$. Then the fractional integral operator \eqref{6.1} yields
\begin{align}
{\bf I}^\mu (x^{A-I}) = \Gamma(A) \Gamma ^{-1}(A+\mu I) x^{A+(\mu - 1)I}.
\end{align}
\end{lemma}
In the next two theorems, we find the fractional order integral and derivative of matrix function ${}_{p}R_{q}(A, B; z)$. 
\begin{theorem}\label{t1}
Let $A$, $B$, $C_1$, $\dots, C_p$, $D_1, \dots, D_q$ be matrices in $\mathbb{C}^{r \times r}$ and $\mu \in \mathbb{C}$ such that $D_i D_j = D_j D_i, 1 \le i,j \le q$ and $\Re(\mu) > 0$. Then the fractional integral of the matrix function ${}_{p}R_{q}(A, B; z)$ is given by
\begin{align}
&{\bf I}^\mu [{}_{p}R_{q}(A, B; z) z^{D_j - I}] \nonumber\\
& = {}_pR_q\left(\begin{array}{c}
C_1, \dots, C_p \\
D_1, \dots, D_{j-1}, D_j + \mu I, D_{j+1}, \dots, D_q
\end{array}\mid A, B; z\right) z^{D_j +(\mu - 1)I} \nonumber\\
& \quad \times \Gamma(D_j) \Gamma^{-1} (D_j + \mu I). 
\end{align}
\end{theorem}
\begin{proof}
From Equation \eqref{6.1}, we have
\begin{align}
&{\bf I}^\mu [{}_{p}R_{q}(A, B; z) z^{D_j - I}]\nonumber\\
 & = \frac{1}{\Gamma(\mu)} \int_{0}^{z} (z-t)^{\mu - 1} {}_{p}R_{q}(A, B; t) t^{D_j - I} dt\nonumber\\
 & = \frac{1}{\Gamma(\mu)} \sum_{n\geq 0}  \, (C_1)_n\dots (C_{p})_n \left(\int_{0}^{z} (z-t)^{\mu - 1} t^{D_j + (n-1) I} dt\right) (D_1)_n^{-1}\dots (D_q)_n^{-1} \, \frac{1}{n!}\nonumber\\
 & = \sum_{n\geq 0}  \, (C_1)_n\dots (C_{p})_n \left({\bf I}^\mu \, z^{D_j + (n-1) I} \right) (D_1)_n^{-1}\dots (D_q)_n^{-1} \, \frac{1}{n!}.
\end{align}
Using the Lemma \ref{l1}, we get
\begin{align}
{\bf I}^\mu [{}_{p}R_{q}(A, B; z) z^{D_j - I}] & = \frac{1}{\Gamma(\mu)} \sum_{n\geq 0}  \, (C_1)_n\dots (C_{p})_n \Gamma(D_j + nI) \Gamma ^{-1}(D_j + nI +\mu I)\nonumber\\
&\quad \times  z^{D_j+(n+\mu - 1)I} (D_1)_n^{-1}\dots (D_q)_n^{-1} \, \frac{1}{n!}\nonumber\\
& = {}_pR_q\left(\begin{array}{c}
C_1, \dots, C_p \\
D_1, \dots, D_{j-1}, D_j + \mu I, D_{j+1}, \dots, D_q
\end{array}\mid A, B; z\right)  \nonumber\\
& \quad \times z^{D_j +(\mu - 1)I} \, \Gamma(D_j) \Gamma^{-1} (D_j + \mu I). 
\end{align}
This completes the proof.
\end{proof}
\begin{theorem}
Let $A$, $B$, $C_1$, $\dots, C_p$, $D_1, \dots, D_q$ be matrices in $\mathbb{C}^{r \times r}$ and $\mu \in \mathbb{C}$ such that $D_i D_j = D_j D_i, 1 \le i,j \le q$ and $\Re(\mu) > 0$. Then the fractional integral of the matrix function ${}_{p}R_{q}(A, B; z)$ is given by
\begin{align}
&{\bf D}^\mu [{}_{p}R_{q}(A, B; z) z^{D_j - I}] \nonumber\\
& = {}_pR_q\left(\begin{array}{c}
C_1, \dots, C_p \\
D_1, \dots, D_{j-1}, D_j - \mu I, D_{j+1}, \dots, D_q
\end{array}\mid A, B; z\right) z^{D_j - (\mu - 1)I} \nonumber\\
& \quad \times \Gamma(D_j) \Gamma^{-1} (D_j - \mu I). \label{6.8}
\end{align}
\end{theorem}
\begin{proof}
	The fractional derivative operator \eqref{6.2} and Theorem \ref{t1} together yield
	\begin{align}
&{\bf D}^\mu [{}_{p}R_{q}(A, B; z) z^{D_j - I}] \nonumber\\
& = \left(\frac{d}{dz}\right)^r {}_pR_q\left(\begin{array}{c}
C_1, \dots, C_p \\
D_1, \dots, D_{j-1}, D_j + (r-\mu) I, D_{j+1}, \dots, D_q
\end{array}\mid A, B; z\right) z^{D_j +(r-\mu - 1)I} \nonumber\\
& \quad \times \Gamma(D_j) \Gamma^{-1} (D_j + (r-\mu) I).
	\end{align}
	Now, proceeding exactly in the same manner as in Theorem \ref{et1}, we get \eqref{6.8}.	
\end{proof}
\section{Special Cases}
The  matrix function ${}_{p}R_{q}(A, B; z)$ reduces to several special matrix functions. These matrix functions are considered as matrix generalizations of respective  classical matrix  functions such as the generalized hypergeometric matrix function, the Gauss hypergeometric matrix function, the confluent hypergeometric matrix function, the matrix M-series, the Wright matrix function and the Mittag-Leffler matrix function and its  generalizations.  We also discuss some matrix polynomials  as particular cases. 

We start with the special case $A = B = I$ and $C_p = I$. The matrix function ${}_{p}R_{q}(A, B; z)$ reduces to 
\begin{align}
{}_pR_q\left(\begin{array}{c}
C_1, \dots, C_{p-1}, I \\
D_1, \dots, D_q
\end{array}\mid I, I; z\right)
& = \sum_{n\geq 0}  \, (C_1)_n\dots (C_{p-1})_n \, (D_1)_n^{-1}\dots (D_q)_n^{-1} \, \frac{z^n}{n!}\nonumber\\
& = {}_{p-1} F_q (C_1, \dots, C_{p-1}, D_1, \dots, D_q; z),
\end{align}
which is known as generalized hypergeometric matrix function with $p-1$ matrix parameters in the numerator and $q$ in the denominator \cite{ds1}. For $C_1 = A_1, C_2 = B_1, C_3 = I, D_1 = C$ and $A = B = I$, the matrix function  ${}_{p}R_{q}(A, B; z)$ reduces to the Gauss hypergeometric matrix function ${}_2F_1(A_1, B_1; C; z)$.
Similarly, for $C_1 = A_1, C_2 = I, D_1 = C$ and $A = B = I$, ${}_{p}R_{q}(A, B; z)$ reduces to the confluent hypergeometric matrix function ${}_1F_1(A_1; C; z)$. 
 
 For $C_p = I$, the matrix function ${}_{p}R_{q}(A, B; z)$ leads to  the matrix analogue of the generalized $M$-series \cite{sj}. 
\begin{align}
{}_pR_q\left(\begin{array}{c}
C_1, \dots, C_{p-1}, I \\
D_1, \dots, D_q
\end{array}\mid A, B; z\right)
& = \sum_{n\geq 0}  \Gamma^{-1}(nA + B)  (C_1)_n\dots (C_{p-1})_n \nonumber\\
& \quad \times  (D_1)_n^{-1}\dots (D_q)_n^{-1} \, z^n\nonumber\\
& = {}_{p-1} M^{(A, B)}_q (C_1, \dots, C_{p-1}, D_1, \dots, D_q; z).\label{4.2}
\end{align}
With $p = 1, q = 0$, $C_1 = I$ and $B = I$, the matrix function ${}_{p}R_{q}(A, B; z)$ reduces to
\begin{align}
{}_1R_0\left(\begin{array}{c}
 I \\
-
\end{array}\mid A, I; z\right)
& = \sum_{n\geq 0}  \Gamma^{-1}(nA + I)  z^n = E_A (z),\label{4.3}
\end{align}
for $p = 1, q = 0$ and $C_1 = I$, the matrix function ${}_{p}R_{q}(A, B; z)$ gives
\begin{align}
{}_1R_0\left(\begin{array}{c}
I \\
-
\end{array}\mid A, B; z\right)
& = \sum_{n\geq 0}  \Gamma^{-1}(nA + B) \,  z^n = E_{A, B} (z),\label{4.4}
\end{align}
with one matrix parameter, $C_1 = C$, ${}_{p}R_{q}(A, B; z)$ becomes
\begin{align}
{}_1R_0\left(\begin{array}{c}
C \\
-
\end{array}\mid A, B; z\right)
& = \sum_{n\geq 0}  \Gamma^{-1}(nA + B) \, (C)_n \,  \frac{z^n}{n!} = E^{C}_{A, B} (z)\label{4.5}
\end{align}
and for two numerator matrix parameter, $C_1 = C, \, C_2 = I$ and one denominator matrix parameter $D_1 = D$, ${}_{p}R_{q}(A, B; z)$ reduces to
\begin{align}
{}_2R_1\left(\begin{array}{c}
C, I \\
D
\end{array}\mid A, B; z\right)
& = \sum_{n\geq 0}  \Gamma^{-1}(nA + B) \, (C)_n \, (D)_n^{-1}  z^n = E^{C, D}_{A, B} (z).\label{4.6}
\end{align}
We define the matrix functions obtained in \eqref{4.3}-\eqref{4.6} as the matrix analogue of the classical Mittag-Leffler function \cite{ml}, Wiman's function \cite{aw}, the generalized Mittag-Leffler function in three parameters \cite{trp} and the generalized Mittag-Leffler function in four parameters \cite{tos}, respectively.

For $p = q = 0$, with replacement of $B$ by $B+I$ and $z$ by $-z$, the matrix function ${}_{p}R_{q}(A, B; z)$ turns into the
generalized Bessel-Maitland matrix function \cite{snd}
\begin{align}
{}_0R_0\left(\begin{array}{c}
- \\
-
\end{array}\mid A, B+I; -z\right)
& = \sum_{n\geq 0} \frac{ \Gamma^{-1}(nA +B+I) \,  (-z)^n}{n!} = J_{A}^{B}(z).  \label{e4.4}
\end{align}

Matrix polynomials such as the Jacobi matrix polynomial, the generalized Konhauser matrix polynomial, the Laguerre matrix polynomial, the Legendre matrix polynomial, the Chebyshev matrix polynomial and the Gegenbauer matrix polynomial can be presented as particular cases of the matrix function  ${}_{p}R_{q}(A, B; z)$. The matrix polynomial dependency chart is given below:
 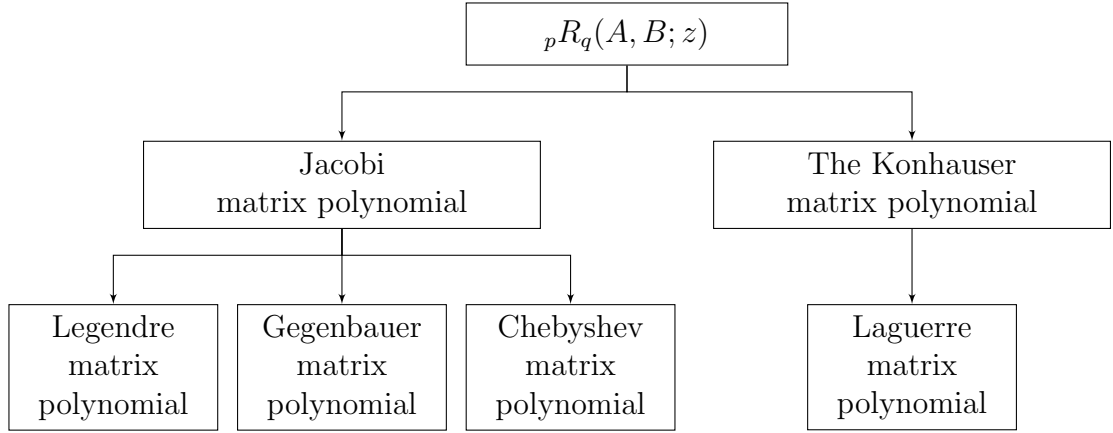
\begin{figure}[!htb]	
 	\caption{Special cases}
 	\vspace*{0.5cm}
 	\begin{tikzpicture}[
 	>=latex',
 	auto
 	]
 	\node [intl] (ki1)  {${}_{p}R_{q}(A, B; z)$};
 	\node [int]  (ki2) [node distance=1cm and -1cm,below right=of ki1] { The  Konhauser \\ matrix  polynomial };
 	\node [intg]  (ki4) [node distance=1 cm and -1cm,below   =of ki2] {Laguerre matrix polynomial };
 	\node [int]  (ki3) [node distance=1cm and -1cm,below left=of ki1] { Jacobi \\ matrix polynomial};
 	\node [intg]  (ki5) [node distance=1cm and -1cm,below  right =of ki3] {Chebyshev \\ matrix \\ polynomial };
 	\node [intg]  (ki6) [node distance=1cm and -1cm,below =of ki3] {  Gegenbauer \\ matrix \\ polynomial};
 	\node [intg]  (ki7) [node distance=1cm and -1cm,below left  =of ki3] { Legendre \\ matrix \\ polynomial };
 	
 	\draw[->] (ki1) -- ($(ki1.south)+(0,-0.35)$) -| (ki2);
 	\draw[->] (ki1) -- ($(ki1.south)+(0,-0.35)$) -| (ki3);
 	\draw[->] (ki2) -- ($(ki2.south)+(0,-0.35)$) -| (ki4);
 	\draw[->] (ki3) -- ($(ki3.south)+(0,-0.35)$) -| (ki5);
 	\draw[->] (ki3) -- ($(ki3.south)+(0,-0.35)$) -| (ki6);
 	\draw[->] (ki3) -- ($(ki3.south)+(0,-0.35)$) -| (ki7);
 	\end{tikzpicture}
 \end{figure}

More explicitly, the Jacobi matrix polynomial can be written in term of the matrix function ${}_{p}R_{q}(A, B; z)$, for $p = 2$, $q = 1$, $C_1 = A + C + (k+1)I$, $C_2 = -kI$, $D_1 = C + I$, $A = 0$, $B = C+I$ and $ z = \frac{1+x}{2}$, as
\begin{align}
P_k^{(A, C)} (x) &= \frac{(-1)^k}{k!} {}_2R_1\left(\begin{array}{c}
A + C + (k+1)I, -kI \\
C+I
\end{array}\mid 0, C+I; \frac{1+x}{2}\right)\nonumber\\
 & \quad \times \Gamma(C+(k+1)I). 
\end{align}
For $p = 2$, $q = 1$, $C_1 = (k+1)I$, $C_2 = -kI$, $D_1 = D$, $A = 0$ and $z = \frac{1-x}{2}$, the matrix function ${}_{p}R_{q}(A, B; z)$ reduces to the Legendre matrix polynomial
\begin{align}
P_k (x,D) &=  {}_2R_1\left(\begin{array}{c}
 (k+1)I, -kI \\
D
\end{array}\mid 0, B; \frac{1-x}{2}\right).
\end{align} 
Similarly, the Gegenbauer matrix polynomial in terms of the matrix function ${}_{p}R_{q}(A, B; z)$ can be expressed as
\begin{align}
C^{D}_k (x) &= \frac{(2D)_k}{k!} {}_2R_1\left(\begin{array}{c}
2D+kI, -kI \\
D + \frac{1}{2} I
\end{array}\mid 0,B; \frac{1-x}{2}\right).
\end{align} 
The Konhauser matrix polynomial in terms of the matrix 
\begin{align}
Z^C_m (x, k) &= \frac{\Gamma(C+(km + 1)I)}{\Gamma(m+1)} {}_1R_0\left(\begin{array}{c}
-mI \\
-
\end{array}\mid kI, C + I; x^k\right).\label{7.12}
\end{align} 
The Laguerre matrix polynomial can be obtained by taking $k = 1$ in  Equation \eqref{7.12}.
 
Note that the properties of these matrix functions and polynomials can be deduced from the corresponding properties of the matrix function ${}_{p}R_{q}(A, B; z)$.

\subsection*{Acknowledgements}
The authors thank the referees for valuable suggestions that led to a better presentation of the paper.


\EditInfo{November 03, 2019}{February 01, 2020}{Karl Dilcher}

\end{paper}